\documentclass[11pt]{article}
\usepackage{mathrsfs}
\usepackage{amsmath, bm, psfrag, epsfig}
\usepackage{latexsym}
\usepackage{amssymb}
\usepackage{amsfonts}
\usepackage{appendix}

 \makeatother
\usepackage{color}
\usepackage{cite}
\usepackage{hyperref}
\hypersetup{colorlinks=true,
            linkcolor=cyan,
            urlcolor=red,
            citecolor=green,
            }
\usepackage{geometry}
\usepackage{epsfig}
\usepackage{epstopdf}
\usepackage{float}

\usepackage{enumerate}
\usepackage{hhline}
\usepackage{multirow}
\usepackage{multicol}
\usepackage{booktabs}
\usepackage{lineno}

\usepackage{graphicx}
\newtheorem{thm}{Theorem}[section]

\newtheorem{lemma}[thm]{Lemma}

\begin{document}
\title{\bf

The existence of even factors based on spectral conditions of graphs}
\author{Jiasheng Li, \, Xiaoyun Lv, \, Shou-Jun Xu
\thanks{ Corresponding author.\\ \phantom{aaa}  Email address: shjxu@lzu.edu.cn }\\
School of Mathematics and Statistics,\\Lanzhou University,
Lanzhou, Gansu, P.R. China}
\date{}
\maketitle

\begin{abstract}
Let $G=(V(G),E(G)) $  be a graph with vertex set $V(G)$ and edge set $E(G)$. 
An even factor of $G$ is a spanning subgraph $F$ such that every vertex in $F$ has a nonzero even degree.
Note that $\delta(G)\geq 2$ is a trivial necessary condition for a graph to have an even factor, where \( \delta(G) \) is the minimum degree of \( G \).   
In this paper, for a connected graph $G$ with minimum degree $\delta$, we establish  a lower bound on the signless
  Laplacian spectral radius of $G$ and an upper bound on the distance spectral radius of $G$ such
  that $G$ contains an even factor. 
\medskip


\medskip

\end{abstract}

\section{Introduction}

\hspace{2em}In this paper, we only consider simple and undirected graphs.
 Let $G$  be a graph with vertex set $V(G)$  and edge set $E(G)$. We denote by $|V(G)|=n$  and $|E(G)|=e(G)$   the order and the size of  $G$, respectively. 
 For any subset $S \subseteq V(G)$, let $G[S]$ denote the subgraph of $G$ induced by $S$, and write $G - S = G[V(G) \setminus S]$. 
 The set of all vertices adjacent to \( v \in V(G) \) is called the neighborhood of \( v \), denoted as \( N_G(v) \).
 We denote by \( d_G(v) = |N_G(v)| \) the degree of \( v \in V(G) \) and by \( \delta(G) \) (or \( \delta \) for short) the minimum degree of \( G \), respectively. 
 Let \( o(G) \) denote the number of odd components in \( G \). For two vertex-disjoint graphs \( G_1 \) and \( G_2 \), we use \( G_1 \cup G_2 \) to denote the disjoint union of \( G_1 \) and \( G_2 \). 
 The join \( G_1 \lor G_2 \) is the graph obtained from \( G_1 \cup G_2 \) by adding all possible edges between  edges between \( V(G_1) \) and \( V(G_2) \). 
 For the graph theoretic notation and terminology not defined here, we refer to \cite{06}.\par

Let \( G \) be a connected graph with vertex set \( V(G) = \{v_1, v_2, \cdots, v_n\} \). 
We denote by \( A(G) = (a_{ij})_{n \times n} \) the adjacency matrix of \( G \), where \( a_{ij} = 1 \) if \( v_i \) and \( v_j \) are adjacent in \( G \), and \( a_{ij} = 0 \) otherwise.
 The signless Laplacian matrix of \( G \) is defined as \( Q(G) = A(G) + D(G) \), where \( D(G) = {\rm diag}\) \((d_G(v_1), d_G(v_2), \cdots, d_G(v_n)) \) is the diagonal matrix of vertex degrees of \( G \). 
 The distance between \( v_i \) and \( v_j \) is denoted by \( d_{ij} \), which is the length of a shortest path from \( v_i \) to \( v_j \). 
 The distance matrix of \( G \) is an \( n \times n \) real symmetric matrix whose \((i, j)\)-entry is \( d_{ij} \), denoted by \( \mathcal{D}(G) \).  
 Particularly, the largest eigenvalue \( \rho_Q(G) \) of \( Q(G) \) is called the signless Laplacian spectral radius of \( G \) and  the largest eigenvalue \( \rho_\mathcal{D}(G) \) of \( \mathcal{D}(G) \) is called the distance spectral radius of \( G \). 
 We refer the reader to \cite{21} and \cite{01} for more definitions of the signless Laplacian spectral radius and the distance spectral radius, respectively.\par

A spanning subgraph $F$ of $G$ is called an even factor of $G$ if $ d_F(v)$ is a nonzero even number for each $ v \in V(G)$. 
Xiong \cite{15} and Lv \cite{10} provided necessary and sufficient conditions for the  \( n \)-times iterated line graph \( L^n(G) \) to have an even factor via branch-bonds. 
Zhang \cite{20} and  Xiong \cite{14} characterized forbidden subgraphs for the existence of even factors in a  graph. 
Steffen, Wolf \cite{12} and Fleischner \cite{03} provided sufficient conditions for critical  graphs to contain an even factor. 
Gyula  \cite{11} established a formula for the maximum cardinality of an even factor in weakly symmetric digraphs. 
Yan and Kano \cite{16} claimed that if \( o(G-S) < |S| \) for all \( S \subseteq V(G) \) with \( |S| \geq 2 \), then \( G \) contains an even factor. 
Recently, Zhou, Bian and Wu \cite{18} presented a sufficient condition based on the size and the adjacency spectral radius for a connected graph to contain an  even factor. 
More research on even factors can be found in \cite{02,04,08,09}.\par

Directly motivated by \cite{16} and  \cite{18}, it is natural and interesting to give some sufficient conditions to ensure that a graph contains an even factor from the perspective of signless Laplacian spectral radius and distance spectral radius. 
Our main results are presented as follows.
\begin{thm}\label{t-1.1}
	Let \( G \) be a connected graph of even order \( n \geq \max\{7\delta - 7, \frac{1}{4}\delta^2 + \frac{1}{2}\delta + 6\} \), where \( \delta \geq 2 \) is the minimum degree of \( G \). If
	\begin{center}
		$\rho_Q(G) \geq \rho_Q(K_\delta \lor (K_{n-2\delta+1} \cup (\delta - 1)K_1)),$
	\end{center}
	then \( G \) contains an even factor, unless \( G \cong K_\delta \lor (K_{n-2\delta+1} \cup (\delta - 1)K_1)\).
	
\end{thm}

\begin{thm}\label{t-1.2}
Let \( G \) be a connected graph of even order \( n \geq \max\{8\delta - 7, \frac{1}{3}\delta^2 + 3\} \),
where \( \delta \geq 2 \) is the minimum degree of \( G \). If
\begin{center}
	$\rho_\mathcal{D}(G) \leq \rho_\mathcal{D}(K_\delta \lor (K_{n-2\delta+1} \cup (\delta - 1)K_1)),$
\end{center}
then \( G \) contains an even factor, unless \( G \cong K_\delta \lor (K_{n-2\delta+1} \cup (\delta - 1)K_1)\).
\end{thm}

\section{Preliminary lemmas}
\hspace{2em}In this section, we list a series of lemmas that will be employed in the subsequent sections. 

\begin{lemma}[\hspace{-0.1pt}\cite{13}\hspace{-0.1pt}]\label{t-2.1}
	Let $H$ be a subgraph of a connected graph $G$, then $\rho_Q(H) \leq \rho_Q(G)$. If
	$H$ is a proper subgraph of a connected graph $G$, then $\rho_Q(H) < \rho_Q(G)$.
\end{lemma}

\begin{lemma}[\hspace{-0.1pt}\cite{5}\hspace{-0.1pt}]\label{t-2.2}
Let $e$ be an edge of a connected graph $G$ such that $G-e$ is also connected, 
then $\rho_\mathcal{D}(G) < \rho_\mathcal{D}(G-e)$.
\end{lemma}

\begin{lemma}[\hspace{-0.1pt}\cite{19}\hspace{-0.1pt}]\label{t-2.3}
Let \hspace{0.1em}$n=\sum_{i=1}^{t} n_i +s$. If $n_1 \geq n_2 \geq \cdots \geq n_t \geq p \geq 1$ and  $n_1 < n - s - p(t - 1)$, then
\[
\rho_Q(K_s \vee (K_{n_1} \cup K_{n_2} \cup \cdots \cup K_{n_t})) < \rho_Q(K_s \vee (K_{n-s-p(t-1)} \cup (t-1)K_p)).
\]
\end{lemma}

\begin{lemma}[\hspace{-0.1pt}\cite{19}\hspace{-0.1pt}]label{t-2.4}
	Let \hspace{0.1em}$n=\sum_{i=1}^{t} n_i +s$. If $n_1 \geq n_2 \geq \cdots \geq n_t \geq p \geq 1$ and  $n_1 < n - s - p(t - 1)$, then  
	\[
	\rho_\mathcal{D}(K_s \vee (K_{n_1} \cup K_{n_2} \cup \cdots \cup K_{n_t})) > \rho_\mathcal{D}(K_s \vee (K_{n-s-p(t-1)} \cup (t-1)K_p)).
	\]
\end{lemma}\par

Let \( M \) be a real symmetric matrix of order $n$ and $V = \{1, 2, \cdots, n\}$.
Given a partition $\Pi: V = V_1 \cup V_2 \cup \cdots \cup V_r$, the matrix $M$ can be correspondingly partitioned as
\[
M = 
\begin{pmatrix}
	M_{11} & \cdots & M_{1r} \\
	\vdots & \ddots & \vdots \\
	M_{r1} & \cdots & M_{rr}
\end{pmatrix},
\] 
where $M_{ij}$ $(i,j = 1,2,\cdots,r)$ denotes the submatrix of $M$ formed by the rows in $V_i$ and the columns in $V_j$.  
Let $q_{ij}$ denote the average row sum of $M_{ij}$. 
Then the matrix $A_{\Pi} = (q_{ij})_{r \times r}$ is called the quotient matrix of $M$ with respect to $\Pi$. 
The partition $\Pi$ is said to be equitable if each block $M_{ij}$ of $M$ has constant row sum $q_{ij}$. 
Also, we say that the quotient matrix $M_\Pi$ is equitable if $\Pi$ is an equitable partition of $M$. 
Particularly, if  \( M \) is nonnegative and irreducible, then by the Perron--Frobenius theorem \cite{07}, there exists a unique unit positive eigenvector, say \( \mathbf{x} = (x_1, x_2, \cdots, x_n)^T \) of \( M \), corresponding to the largest eigenvalue  of \( M \). 
As usual, \( \mathbf{x} \) is said to be the Perron vector of \( M \).

\begin{lemma}[\hspace{-0.1pt}\cite{17}\hspace{-0.1pt}]\label{t-2.5}
	Let \( M \) be a real symmetric matrix, and let \(\lambda_1(M)\) be the largest eigenvalue of \( M \). If \( A_{\Pi} \) is an equitable quotient matrix of \( M \), then the eigenvalues of \( A_{\Pi} \) are also eigenvalues of \( M \). Furthermore, if \( M \) is nonnegative and irreducible, then \(\lambda_1(M) = \lambda_1(A_{\Pi})\).
\end{lemma}\par

Note that the Wiener index of a connected graph \( G \) with \( n \) vertices is denoted by \( W(G) = \sum_{i<j} d_{ij} \). 
The following lemma can be easily obtained by the Rayleigh quotient \cite{07}.

\begin{lemma}[\hspace{-0.1pt}\cite{21}\hspace{-0.1pt}]\label{t-2.6}
Let \( G \) be a connected graph of order \( n \). Then
\[
\rho_\mathcal{D}(G) = \max_{\mathbf{x} \neq \mathbf{0}} \frac{\mathbf{x}^T \mathcal{D}(G)\mathbf{x}}{\mathbf{x}^T \mathbf{x}} \geq \frac{\mathbf{1}^T \mathcal{D}(G)\mathbf{1}}{\mathbf{1}^T \mathbf{1}} = \frac{2W(G)}{n},
\]
where \( \mathbf{1} = (1, 1, \cdots, 1)^T \).
\end{lemma}

\begin{lemma}\label{t-2.7}
	Let \( n, \delta \) be two positive integers with \( n \geq 7\delta - 7 \) and \( \delta \geq 2 \). 
	Then
	\[
	2n - 2\delta < \rho_Q(K_\delta \vee (K_{n-2\delta+1} \cup (\delta - 1)K_1)) < 2n - \delta.
	\]
\end{lemma}
\textbf{Proof.} Let $G_{*}=K_{\delta}\vee(K_{n-2\delta+1}\cup(\delta-1)K_{1})$. 
It is clear that $K_{n-\delta+1}$ is a proper subgraph of $G_{*}$, this implies $\rho_{Q}(G_{*})>2n-2\delta$ by Lemma 2.1. 
We partition $V(G_{*})$ into $V(K_{\delta})\cup V(K_{n-2\delta+1})\cup V((\delta-1)K_{1})$. 
According to this partition, the quotient matrix of $Q(G_{*})$ is
\[
A^{*}_{\Pi}=\begin{bmatrix}
	n+\delta-2 & n-2\delta+1 & \delta-1 \\ 
	\delta & 2n-3\delta & 0 \\ 
	\delta & 0 & \delta
\end{bmatrix}.
\]
By a simple calculation, the characteristic polynomial of $A^{*}_{\Pi}$ is
\[
\varphi(A^{*}_{\Pi},x)= x^{3}+(-3n+\delta+2)x^{2}+(2n^{2}+\delta n-4n-4\delta^{2}+4\delta)x -2\delta n^{2}+4\delta^{2}n+2\delta n-2\delta^{3}-2\delta^{2}.
\]
We will prove that $\varphi(A^{*}_{\Pi},x)$ is increasing on the interval $[2n-\delta,+\infty)$. 
Since
\begin{align*}
	\varphi'\left(A^{*}_{\Pi},x\right) &= 3x^{2}+(2\delta+4-6n)x-4\delta^{2}+\delta n+4\delta+2n^{2}-4n \\
	&\geq 2n^{2}+(4-\delta)n-3\delta^{2} \quad (\text{since } x\geq 2n-\delta) \\
	&\geq (88\delta-161)\delta+70 \quad (\text{since } n\geq 7\delta-7) \\
	&> 0 \quad (\text{since } \delta\geq 2),
\end{align*}
 it follows that
\begin{align*}
	\varphi(A^{*}_{\Pi},x) &\geq \varphi(A^{*}_{\Pi},2n-\delta) \\
	&= 2\delta n^{2}+(6\delta-6\delta^{2})n+2\delta^{3}-4\delta^{2} \\
	&\geq (58\delta-116)\delta^{2}+56\delta \quad (\text{since } n\geq 7\delta-7) \\
	&> 0 \quad (\text{since } \delta\geq 2).
\end{align*}
Notice that the partition $V(G_{*})=V(K_{\delta})\cup V(K_{n-2\delta+1})\cup V((\delta-1)K_{1})$ is equitable, we conclude that $\rho_{Q}({G_{*}})$ equals the largest root of $\varphi(A^{*}_{\Pi},x)=0$  by Lemma 2.5. 
Concluding above results, we deduce that $\rho_{Q}(K_{\delta}\vee(K_{n-2\delta+1}\cup(\delta-1)K_{1}))<2n-\delta$, as required.

\begin{lemma}[\cite{11}]\label{t-2.8}
	Let \( G \) be a graph of even order \( n \). 
	Then \( G \) contains an even factor if
	\[
	o(G - S) < |S|
	\]
   for all \( S \subseteq V(G) \) with \( |S| \geq 2 \), where \( o(G - S) \) denotes the number of odd components of \( G - S \).
\end{lemma}

\section{Proof of Theorem 1.1}
\hspace{2em}In this section, we prove Theorem 1.1.\\
\textbf{Proof.} Suppose to the contrary that $G$ has no even factor. 
By virtue of Lemma 2.8, we have 
\[ 
o(G-S) \geq |S| \tag{1}
\]for some subset $S \subseteq V(G)$ with $ |S| \geq 2$. For convenience, let $|S|=s$.
Choose a connected graph $G$ with even order $n$ such that its signless Laplacian spectral radius is as large as possible. 
It is clear that $G$ is a spanning subgraph of $G_1 = K_s \vee (K_{n_1} \cup K_{n_2} \cup \cdots \cup K_{n_s})$ for some odd integers $n_1 \geq n_2 \geq \cdots \geq n_s$  with  $\sum_{i=1}^s n_i = n - s.$  
By Lemma 2.1, we can deduce that 
\[
\rho_Q(G) \leq \rho_Q(G_1), \tag{2}
\] with equality occurring if and only if $G\cong G_{1}$. 
We divide the proof into the following four cases.\\
\textbf{Case 1.} $s\geq\delta+1$\par
Let $G_{2}=K_{s}\vee(K_{n-2s+1}\cup(s-1)K_{1})$, where $n\geq 2s$. 
According to  Lemma 2.3, we conclude that
\[
\rho_{Q}(G_{1})\leq\rho_{Q}(G_{2}), \tag{3}
\] with equality holding if and only if  $(n_{1},n_{2},\cdots,n_{s})=(n-2s+1,1,\cdots,1)$. 
 We partition the vertex set of  $G_{2}=K_{s}\vee(K_{n-2s+1}\cup(s-1)K_{1})$ into $V(K_{s})\cup V(K_{n-2s+1})\cup V((s-1)K_{1})$. 
 The quotient matrix corresponding to the partition is  
\[
A^{2}_{\Pi}=\begin{bmatrix}
	n+s-2 & n-2s+1 & s-1 \\ 
	s & 2n-3s & 0 \\ 
	s & 0 & s
\end{bmatrix}.
\]
By a simple calculation, the characteristic polynomial of $A^{2}_{\Pi}$ is
\begin{align*}
	\varphi(A^{2}_{\Pi},x) &= x^{3}+(-3n+s+2)x^{2}+(2n^{2}+sn-4n-4s^{2}+4s)x \\
	&\quad -2sn^{2}+4s^{2}n+2sn-2s^{3}-2s^{2}.
\end{align*}
Note that the partition $V(G_{2})=V(K_{s})\cup V(K_{n-2s+1})\cup V((s-1)K_{1})$ is equitable. 
According to Lemma 2.5, the largest root of $\varphi(A^{2}_{\Pi},x)=0$ equals $\rho_{Q}(G_{2})$.\par
Recall that $G_{*}=K_{\delta}\vee(K_{n-2\delta+1}\cup(\delta-1)K_{1})$, and the characteristic polynomial of $A^{*}_{\Pi}$ is
\begin{align*}
	\varphi(A^{*}_{\Pi},x) &= x^{3}+(-3n+\delta+2)x^{2}+(2n^{2}+\delta n-4n-4\delta^{2}+4\delta)x \\
	&\quad -2\delta n^{2}+4\delta^{2}n+2\delta n-2\delta^{3}-2\delta^{2}.
\end{align*}
Observe that the partition $V(G_{*})=V(K_{\delta})\cup V(K_{n-2\delta+1})\cup V((\delta-1)K_{1})$ is equitable. 
According to Lemma 2.5, the largest root of $\varphi(A^{*}_{\Pi},x)=0$ equals $\rho_{Q}(G_{*})$.\\
By a simple calculation, we obtain
\[
\varphi(A^{2}_{\Pi},x)-\varphi(A^{*}_{\Pi},x) 
= (s-\delta)\eta_{1}(x), \tag{4}
\]
where
\begin{center}
	$\eta_{1}(x) = x^{2}+(n+4-4s-4\delta)x
	-2n^{2}+2n+4sn+4\delta n-2s^{2}-2\delta^{2}-2s-2\delta-2s\delta.$
\end{center}
Notice that the symmetry axis of $\eta_{1}(x)$ satisfies
\[
\frac{n+4-4s-4\delta}{-2} \leq 2n-2\delta 
\]
by $s\leq\frac{n}{2}$ and $\delta\geq 2$, this implies that $\eta_{1}(x)$ is increasing with respect to $x \geq 2n-2\delta$. 
It follows that
\begin{align*}
	\eta_{1}(x) &\geq \eta_{1}(2n-2\delta) \\
	&= -2s^{2}-(4n-6\delta+2)s+4n^{2}-14\delta n+10n+10\delta^{2}-10\delta \\
	&\geq \frac{3}{2}n^{2}+(9-11\delta)n+10\delta^{2}-10\delta \quad (\text{since } s\leq\frac{n}{2}) \\
	&\geq \frac{13}{2}\delta^{2}-17\delta+\frac{21}{2} \quad (\text{since } n\geq 7\delta-7) \\
	&\geq \frac{5}{2} \quad (\text{since } \delta\geq 2) \\
	&> 0. \tag{5}
\end{align*}
From (4) and (5), we conclude that $\varphi(A^{2}_{\Pi},x)>\varphi(A^{*}_{\Pi},x)$ for $x \geq 2n-2\delta$. 
Recall that $\rho_{Q}(G_{*})>2n-2\delta$ by Lemma 2.7, this implies $\rho_{Q}(G_{2})<\rho_{Q}(G_{*})$. 
Together with (2) and (3), we deduce that
\begin{center}
	$\rho_{Q}(G) \leq \rho_{Q}(G_{1}) 
	\leq \rho_{Q}(G_{2}) 
	<\rho_{Q}(G_{*})= \rho_{Q}(K_{\delta} \vee (K_{n-2\delta+1} \cup (\delta-1)K_{1})),$
\end{center}
which contradicts with $\rho_{Q}(G) \geq \rho_{Q}(K_{\delta} \vee (K_{n-2\delta+1} \cup (\delta-1)K_{1}))$.\\
\textbf{Case 2.} $s=\delta$\par
In this case, we have
\begin{center}
	$	\rho_{Q}(G) \leq \rho_{Q}(G_{1}) 
	\leq \rho_{Q}(G_{2})=\rho_{Q}(G_{*})
	= \rho_{Q}(K_{\delta} \vee (K_{n-2\delta+1} \cup (\delta-1)K_{1})),$
\end{center}
with equality holding if and only if $G \cong K_{\delta} \vee (K_{n-2\delta+1} \cup (\delta-1)K_{1})$, a contradiction.\\
\textbf{Case 3.} $3\leq s\leq\delta-1$\par
Let $G_{3}=K_{s} \vee (K_{n-s-(\delta+1-s)(s-1)} \cup (s-1)K_{\delta+1-s})$. 
Recall that $G$ is a spanning subgraph of $G_{1}=K_{s}\vee(K_{n_{1}}\cup K_{n_{2}}\cup\cdots\cup K_{n_{s}})$, where $ n_{1}\geq n_{2}\geq\cdots\geq n_{s}$ are some odd integers with  $\sum_{i=1}^{s}n_{i}=n-s$.
Since $\delta(G_{1})\geq \delta(G)=\delta$, we have $n_{s}\geq\delta+1-s$. 
According to Lemma 2.1 and Lemma 2.3, we deduce that
\[
\rho_{Q}(G) \leq \rho_{Q}(G_{1}) 
\leq \rho_{Q}(G_{3}), \tag{6}
\]
 where the second equality holds if and only if $(n_{1},n_{2},\cdots,n_{s})=(n-s-(\delta+1-s)(s-1),\delta+1-s,\cdots,\delta+1-s)$.\par
We partition $V(G_{3})$ into $V(K_{s}) \cup V(K_{n-s-(\delta+1-s)(s-1)}) \cup V((s-1)K_{\delta+1-s})$. 
Then $Q(G_{3})$ has the equitable quotient matrix
\[
A^{3}_{\Pi}=\begin{bmatrix}
	n+s-2 & n-s-(\delta+1-s)(s-1) & (s-1)(\delta+1-s) \\ 
	s & 2n-s-2(\delta+1-s)(s-1)-2 & 0 \\ 
	s & 0 & 2\delta-s
\end{bmatrix},
\]
and its characteristic polynomial is
\begin{align*}
	\varphi(A^{3}_{\Pi},x) &= x^{3}+(-3n-2s^{2}+2\delta s+5s-4\delta+2)x^{2}+(-4\delta^{2}s+4\delta^{2}-2\delta sn+8\delta n+4\delta s^{2}\\
	&\quad -4\delta s-8\delta+2n^{2}+2ns^{2}-7ns-4n-4s^{2}+12s)x+4s\delta^{2}n-4\delta^{3}n-2\delta^{4}s^{3}\\
	&\quad
	+6\delta^{3}s^{2}-12\delta^{3}s+8\delta^{3}-4\delta n^{2}-4\delta s^{2}n+8s\delta n+8\delta n +4\delta s^{4}-16\delta s^{3}+28\delta s^{2}\\
	&\quad -24\delta s+2sn^{2}-6ns-2s^{3}+10s^{4}-18s^{3}+14s^{2}.
\end{align*}\par
Recall that $G_{*}=K_{\delta} \vee (K_{n-2\delta+1} \cup (\delta-1)K_{1})$, and the characteristic polynomial of $A^{*}_{\Pi}$ is
\begin{align*}
	\varphi(A^{*}_{\Pi},x) &= x^{3}+(-3n+\delta+2)x^{2}+(2n^{2}+\delta n-4n-4\delta^{2}+4\delta)x \\
	&\quad -2\delta n^{2}+4\delta^{2}n+2\delta n-2\delta^{3}-2\delta^{2}.
\end{align*}
By a simple calculation, we get
\[
\varphi(A_{\Pi}^{3},x)-\varphi(A_{\Pi}^{*},x) 
= (\delta-s)\eta_{2}(x), \tag{7}
\]
where
\begin{align*}
	\eta_{2}(x) &= (2s-5)x^{2}+(7n-12-4\delta s-2ns+8\delta+4s)x \\
	&\quad +2\delta^{2}+4\delta sn-8\delta n-2n^{2}+6n+2s^{4}-2\delta s^{3}+6\delta s^{2} \\
	&\quad -10s^{3}+18s^{2}-10\delta s-14s+10\delta. \tag{8}
\end{align*}
Notice that the symmetry axis of $\eta_{2}(x)$ satisfies
\[
\frac{7n-12-4\delta s-2ns+8\delta+4s}{-2(2s-5)} \leq 2n-2\delta 
\]
by $3\leq s\leq\delta-1$ and $n\geq\frac{1}{4}\delta^{2}+\frac{1}{2}\delta+6$, this implies $\eta_{2}(x)$ is increasing with respect to $x \geq 2n-2\delta$, and hence
\begin{align*}
	\eta_{2}(x) &\geq \eta_{2}(2n-2\delta) \\
	&= (4s-8)n^{2}+(34\delta+8s-16\delta s-18)n+2s^{4}-10s^{3}-2\delta s^{3} \\
	&\quad +6\delta s^{2}+18s^{2}+16\delta^{2}s-18\delta s-14s-36\delta^{2}+34\delta.
\end{align*}
Let
\begin{align*}
	f(n) &= (4s-8)n^{2}+(34\delta+8s-16\delta s-18)n+2s^{4}-10s^{3}-2\delta s^{3} \\
	&\quad +6\delta s^{2}+18s^{2}+16\delta^{2}s-18\delta s-14s-36\delta^{2}+34\delta. \tag{9}
\end{align*}
Then we can easily get that
\[
f'(n) = (8s-16)n+34\delta+8s-16\delta s-18.
\]
For $n\geq 2\delta-1$, we have
$f'(n) \geq f'(2\delta-1) 
= 2\delta-2 
> 0.$ 
That is, $f(n)$ is increasing on the interval $[2\delta-1,+\infty)$. 
Since $n\geq\frac{1}{4}\delta^{2}+\frac{1}{2}\delta+6>2\delta-1$, 
we obtain 
\begin{align*}
	f(n) &\geq f\left(\frac{1}{4}\delta^{2}+\frac{1}{2}\delta+6\right) \\
	&= \left(\frac{s}{4}-\frac{1}{2}\right)\delta^{4}+\left(\frac{13}{2}-3s\right)\delta^{3}+(23s-\frac{95}{2})\delta^{2} \\
	&\quad +(-2s^{3}+6s^{2}-86s+181)\delta+2s^{4}-10s^{3}+18s^{2}+178s-396 \\
	&\geq \frac{1}{4}s^{5}-\frac{5}{2}s^{4}+14s^{3}-55s^{2}+\frac{863}{4}s-\frac{513}{2} \quad (\text{since } \delta\geq s+1) \\
	&\geq 132 \quad (\text{since } s\geq 3) \\
	&> 0. \tag{10}
\end{align*}
From (7) $\sim$ (10) and $\delta\geq s+1$, we deduce that $ \varphi(A_{\Pi}^{3},x)>\varphi(A_{\Pi}^{*},x) $ for $x \geq 2n-2\delta$.
Recall that $\rho_{Q}(G_{*})>2n-2\delta$ by Lemma 2.7, we conclude that $\rho_{Q}(G_{3})<\rho_{Q}(G_{*})$. 
It follows that
\begin{center}
	$\rho_{Q}(G) \leq \rho_{Q}(G_{1}) 
	\leq \rho_{Q}(G_{3}) 
	< \rho_{Q}(G_{*}) = \rho_{Q}(K_{\delta} \vee (K_{n-2\delta+1} \cup (\delta-1)K_{1})),$
\end{center}
which contradicts with $\rho_{Q}(G) \geq \rho_{Q}(K_{\delta} \vee (K_{n-2\delta+1} \cup (\delta-1)K_{1}))$.\\
\textbf{Case 4.} $s = 2$\par
In this case, by plugging the value $s=2$ into (8), we get $\eta_{2}(x) = -x^{2}+(3n-4)x+2\delta^{2}-2\delta-2n^{2}+6n-4$. 
Note that $\frac{3n-4}{2} \leq 2n-2\delta  < 2n-\delta$ by $n \geq 7\delta-7$, this implies $\eta_{2}(x)$ is decreasing on the interval $[2n-2\delta,2n-\delta]$. 
It follows that
\begin{align*}
	\eta_{2}(x) &\geq \eta_{2}(2n-\delta) \\
	&= (\delta-2)n+\delta^{2}+2\delta-4 \\
	&\geq 8\delta^{2}-19\delta+10 \quad (\text{since } n \geq 7\delta-7) \\
	&\geq 25 \quad (\text{since } \delta \geq 3) \\
	&> 0.\tag{11}
\end{align*}
For $\delta = 2$, we can directly get that $\eta_{2}(x) \geq \eta_{2}(2n-\delta) = (\delta-2)n+\delta^{2}+2\delta-4 = 4 > 0$ from (11).
Therefore, $\varphi(A_{\Pi}^{3},x)>\varphi(A_{\Pi}^{*},x)$ on the interval $[2n-2\delta,2n-\delta]$.\par
By plugging the value \( s = 2 \) into \( \varphi(B_{\Pi}^3, x) \) in case 3, we have
\begin{align*}
\varphi(B_{\Pi}^3, x) &= x^3 + (4 - 3n)x^2 + (2n^2 + 4\delta n - 10n - 4\delta^2 + 8)x \\
&\quad + 4n^2 - 4\delta n^2 + 4\delta^2 n + 8\delta n - 12n - 8\delta^2 + 8.
\end{align*}
We take the derivative of \( \varphi(B_{\Pi}^3, x) \). If \( x \geq 2n - \delta \), then we obtain that
\[
\begin{aligned}
\varphi'(B_{\Pi}^3, x) &= 3x^2 + (8 - 6n)x + 2n^2 + 4\delta n - 10n - 4\delta^2 + 8 \\
&\geq 2n^2 + (6 - 2\delta)n - \delta^2 - 8\delta + 8 \quad (\text{since } x \geq 2n - \delta) \\
&\geq (83\delta - 148)\delta + 64 \quad (\text{since } n \geq 7\delta - 7) \\
&> 0 \quad (\text{since } \delta \geq 2).
\end{aligned}
\]
It follows that \( \varphi(B_{\Pi}^3, x) \) is increasing with respect to \( x \geq 2n - \delta \), hence
\[
\begin{aligned}
\varphi(B_{\Pi}^3, x) &\geq \varphi(B_{\Pi}^3, 2n - \delta) \\
&= 2\delta n^2 + (-5\delta^2 + 2\delta + 4)n + 3\delta^3 - 4\delta^2 - 8\delta + 8 \\
&\geq 66\delta^3 - 151\delta^2 + 104\delta - 20 \quad (\text{since } n \geq 7\delta - 7) \\
&\geq 112 \quad (\text{since } \delta \geq 2) \\
&> 0.
\end{aligned}
\]
According to Lemma 2.7, we have $2n-2\delta<\rho_{Q}(G_{*})<2n-\delta$. 
From the above discussion, we conclude that  $\rho_{Q}(G_{3})<\rho_{Q}(G_{*})$. 
Togehter with (6), we have
\begin{center}
	$\rho_{Q}(G) \leq \rho_{Q}(G_{1}) 
	\leq \rho_{Q}(G_{3}) 
	<\rho_{Q}(G_{*})= \rho_{Q}(K_{\delta} \vee (K_{n-2\delta+1} \cup (\delta-1)K_{1})),$
\end{center}
a contradiction.\par
This completes the proof of Theorem 1.1.

\section{Proof of Theorem 1.2}
\hspace{2em}In this section, we prove Theorem 1.2.\\
\textbf{Proof.} Suppose to the contrary that $G$ has no even factor. According to Lemma 2.8, we have  
\[
o(G - S) \geq |S| \tag{12}
\]  
for some subset $S \subseteq V(G)$ with $s = |S| \geq 2$. 
Choose a connected graph $G$ with even order $n$ such that its distance spectral radius is as small as possible. 
Recall that $G$ is a spanning subgraph of  $G_1 = K_s \vee (K_{n_1} \cup K_{n_2} \cup \cdots \cup K_{n_s})$, 
where $n_1 \geq n_2 \geq \cdots \geq n_s$  are some odd integers with $\sum_{i=1}^s n_i = n - s.$
By Lemma 2.2, we can deduce that 
\[
\rho_\mathcal{D}(G) \geq \rho_\mathcal{D}(G_1), \tag{13}
\]  
with equality occurring if and only if $G \cong G_1$. 
We consider the following four cases.\\
\textbf{Case 1.} $s \geq \delta + 1$\par
Recall that $G_2 = K_s \vee (K_{n-2s+1} \cup (s-1)K_1)$, where $n \geq 2s$. 
By Lemma 2.4, we obtain that
\[
\rho_\mathcal{D}(G_1) \geq \rho_\mathcal{D}(G_2), \tag{14}
\]
with equality holding if and only if $(n_{1},n_{2},\cdots,n_{s})=(n-2s+1,1,\cdots,1)$. 
We partition $V(K_s \vee (K_{n-2s+1} \cup (s-1)K_1))$ into $V(K_s) \cup V(K_{n-2s+1}) \cup V((s-1)K_1)$. 
Then $\mathcal{D}(G_2)$ has the  quotient matrix
\[
B_{\Pi}^2 = 
\begin{bmatrix}
	n - 2s & s & 2(s-1) \\
	n - 2s + 1 & s - 1 & s - 1 \\
	2(n - 2s + 1) & s & 2(s-2)
\end{bmatrix} .
\]
By a simple calculation, the characteristic polynomial of $B_{\Pi}^2$ is
\[
\varphi(B_{\Pi}^2, x) = x^3 + (5 - n - s)x^2 + (5s^2 - n - 2sn - 8s + 8)x + s^2n - 3sn - 8s + 8s^2 - 2s^3 + 4 .
\]
Note that the partition $V(G_2) = V(K_s) \cup V(K_{n-2s+1}) \cup V((s-1)K_1)$ is equitable, it follows from 
 Lemma 2.5 that $\rho_\mathcal{D}(G_2)$ is the largest root of $\varphi(B_{\Pi}^2,x)=0$.\par
Recall that $G_* = K_\delta \vee (K_{n-2\delta+1} \cup (\delta - 1)K_1)$. 
In terms of the partition $V(G_*)=V(K_\delta) \cup V(K_{n-2\delta+1}) \cup V((\delta - 1)K_1)$, the quotient matrix of $\mathcal{D}(G_{*})$ is equal to
\[
B_{\Pi}^*=\begin{bmatrix}
	n-2\delta & \delta & 2(\delta-1) \\ 
	n-2\delta+1 & \delta-1 & \delta-1 \\ 
	2(n-2\delta+1) & \delta & 2(\delta-2)
\end{bmatrix}
\]
and the characteristic polynomial of $B_{\Pi}^*$ is
\[
\varphi(B_{\Pi}^*,x)= x^{3}+(5-n-\delta)x^{2}+(5\delta^{2}-n-2\delta n-8\delta+8)x + \delta^{2}n-3\delta n-8\delta+8\delta^{2}-2\delta^{3}+4. 
\]
Since the partition $V(G_{*})=V(K_{\delta})\cup V(K_{n-2\delta+1})\cup V((\delta-1)K_{1})$ is equitable. 
According to Lemma 2.5, the largest root of $\varphi(B_{\Pi}^*,x)=0$ equals $\rho_{\mathcal{D}}(G_{*})$.   \par
By a simple calculation, we get
\[
\varphi(B_{\Pi}^2,x)-\varphi(B_{\Pi}^*,x)=(\delta-s)\varphi_{1}(x), \tag{15}
\]
where $\varphi_{1}(x)=x^{2}+(2n+8-5s-5\delta)x+3n+8-sn-\delta n-8s-8\delta+2s^{2}+2s\delta+2\delta^{2}$.\\
According to Lemma 2.6 and $n\geq 8\delta-7$, we have
\[
\rho_{\mathcal{D}}(G_{*})\geq\frac{2W(G_{*})}{n}=\frac{n^{2}+(2\delta-3)n-3\delta^{2}+3\delta}{n}\geq n+\delta-3. \tag{16}
\]
Since $n\geq 2s$, we obtain
\[
\frac{5s+5\delta-2n-8}{2}\leq\frac{n+10\delta-16}{4}< n+\delta-3. \tag{17}
\]
Notice that the symmetry axis of $\varphi_{1}(x)$ is
\[
x=\frac{5s+5\delta-2n-8}{2}. \tag{18}
\]
Combining this with (17), we conclude that $\varphi_{1}(x)$ is increasing on the interval $[n+\delta-3,+\infty)$.
It follows that 
\begin{align*}
	\varphi_1(x) &\geq \varphi_1(n + \delta - 3) \\
	&= 2s^2 + (7 - 6n - 3\delta)s + 3n^2 - 2\delta n - n - 2\delta^2 + 9\delta - 7 \\
	&\geq \frac{1}{2}n^2 + \left(\frac{5}{2} - \frac{7}{2}\delta\right)n - 2\delta^2 + 9\delta - 7 \quad (\text{since } s \leq \frac{n}{2}) \\
	&\geq 2\delta^2 - \frac{5}{2}\delta \quad (\text{since } n \geq 8\delta - 7) \\
	&> 0 \quad (\text{since } \delta \geq 2).
\end{align*}
Combining this with (15) and $s \geq \delta+1$, we get $\varphi(B_{\Pi}^2,x) < \varphi(B_{\Pi}^*,x)$ for $x \geq n+\delta-3$. 
Since $\rho_{\mathcal{D}}(G_{*})\geq n+\delta-3$ by (16), we can deduce that  $\rho_\mathcal{D}(G_2)>\rho_\mathcal{D}(G_*)$.   
From the above disussion, we posses 
\begin{center}
$	\rho_\mathcal{D}(G) \geq \rho_\mathcal{D}(G_1) 
	\geq \rho_\mathcal{D}(G_2) 
	> \rho_\mathcal{D}(G_*) = \rho_\mathcal{D}(K_\delta \vee (K_{n-2\delta+1} \cup (\delta - 1)K_1)),
$
\end{center}
which contradicts with $\rho_\mathcal{D}(G) \leq \rho_\mathcal{D}(K_\delta \vee (K_{n-2\delta+1} \cup (\delta - 1)K_1))$.\\
\textbf{Case 2.} $s = \delta$\par 
In this case, we have 
\begin{center} 
$\rho_\mathcal{D}(G) \geq \rho_\mathcal{D}(G_1) 
	\geq \rho_\mathcal{D}(G_2) =\rho_\mathcal{D}(G_*)
	= \rho_\mathcal{D}(K_\delta \vee (K_{n-2\delta+1} \cup (\delta - 1)K_1))
,$
\end{center}
with equality holding if and only if $G \cong K_\delta \vee (K_{n-2\delta+1} \cup (\delta - 1)K_1)$, a contradiction.\\
\textbf{Case 3.} $3\leq s\leq\delta-1$\par
Recall that $G_{3}=K_{s}\vee(K_{n-s-(\delta+1-s)(s-1)}\cup(s-1)K_{\delta+1-s})$  and $G$ is a spanning subgraph of
$
G_{1}=K_{s}\vee(K_{n_{1}}\cup K_{n_{2}}\cup\cdots\cup K_{n_{s}})$. 
According to Lemma 2.4, we conclude
\begin{align*}
	\rho_{\mathcal{D}}(G) \geq \rho_{\mathcal{D}}(G_{1}) 
	\geq \rho_{\mathcal{D}}(G_{3}), \tag{19}
\end{align*}
where the second equality holds if and only if $(n_{1},n_{2},\cdots,n_{s})=(n-s-(\delta+1-s)(s-1),\delta+1-s,\cdots,\delta+1-s)$.\par
We partition $V(G_{3})$ into $V(K_{s})\cup V(K_{n-s-(\delta+1-s)(s-1)})\cup V((s-1)K_{\delta+1-s})$. Then $\mathcal{D}(G_{3})$ has the equitable quotient matrix
\[
B^{3}_{\Pi}=\begin{bmatrix}
	n-s-(\delta+1-s)(s-1)-1 & s & 2(s-1)(\delta+1-s) \\ 
	n-s-(\delta+1-s)(s-1) & s-1 & (s-1)(\delta+1-s) \\ 
	2n-2s-2(\delta+1-s)(s-1) & s & \delta-s+2(s-2)(\delta+1-s)
\end{bmatrix}
\]
and the characteristic polynomial of $B^{3}_{\Pi}$ is
\begin{align*}
	\varphi(B^{3}_{\Pi},x) &= x^{3}+(2\delta-n-3s-\delta s+s^{2}+5)x^{2}+(6\delta-n-14s+\delta n-13\delta s-3ns+13\delta s^{2}\\
	&\quad -3\delta^{2}s-4\delta s^{3}+2ns^{2}+\delta^{2}+17s^{2}-10s^{3}+2s^{4}+2\delta^{2}s^{2}-2\delta ns+8)x+4\delta-12s\\
	&\quad +\delta n-14\delta s-4ns+20\delta s^{2}-4\delta^{2}s-11\delta s^{3}+2\delta s^{4}+4ns^{2}-ns^{3}+\delta^{2}+21s^{2}-18s^{3} \\
	&\quad +7s^{4}-s^{5}+4\delta^{2}s^{2}-\delta^{3}s^{3}+\delta ns^{2}-3\delta ns+4.
\end{align*}
By a simple calculation, we get
\begin{align*}
	\varphi(B^{3}_{\Pi},x)-\varphi(B_{\Pi}^*,x) 
	= (s-\delta)\varphi_{2}(x), \tag{20}
\end{align*}
where
\begin{align*}
	\varphi_{2}(x) &= (s-3)x^{2} + (2s^{3}-2\delta s^{2}-10s^{2}+2ns+17s+3\delta s+4\delta-14-3n)x\\
	&\quad -s^{4}+7s^{3}+\delta s^{3}-18s^{2}-ns^{2}-4\delta s^{2}+21s+4ns+2\delta s-12-4n+\delta n-2\delta^{2}+7\delta.
\end{align*}\par
We first assume that $s \geq 4$.
Since
\[
\frac{2s^{3}-2\delta s^{2}-10s^{2}+2ns+17s+3\delta s+4\delta-14-3n}{-2(s-3)} \leq n+\delta-3 
\]
by $4\leq s\leq\delta-1$ and $n\geq\frac{1}{3}\delta^{2}+3$, we conclude that $\varphi_{2}(x)$ is increasing on the interval $[n+\delta-3,+\infty)$.
It follows that
\begin{align*}
	\varphi_{2}(x) &\geq \varphi_{2}(n+\delta-3) \\
	&= (3s-6)n^{2}+(2s^{3}-11s^{2}-2\delta s^{2}+7\delta s+9s-4\delta+9)n \\
	&\quad -s^{4}+3\delta s^{3}+s^{3}-2\delta^{2}s^{2}-8\delta s^{2}+12s^{2}+4\delta^{2}s+4\delta s+21s-\delta^{2}-\delta+3 .\tag{21}
\end{align*}
Let
\begin{align*}
	g(n) &= (3s-6)n^{2}+(2s^{3}-11s^{2}-2\delta s^{2}+7\delta s+9s-4\delta+9)n \\
	&\quad -s^{4}+3\delta s^{3}+s^{3}-2\delta^{2}s^{2}-8\delta s^{2}+12s^{2}+4\delta^{2}s+4\delta s+21s-\delta^{2}-\delta+3.\tag{22}
\end{align*}
Then  we can directly get that
\[
g'(n) = (6s-12)n+2s^{3}-11s^{2}-2\delta s^{2}+7\delta s+9s-4\delta+9.
\]
For $n\geq\frac{2\delta s-3\delta+3}{6}$, we have
\begin{align*}
	g'(n) &\geq g'\left(\frac{2\delta s-3\delta+3}{6}\right) \\
	&= 2s^{3}-11s^{2}+12s+2\delta+3 \\
	&\geq 2s^{3}-11s^{2}+14s+5 \quad (\text{since } \delta\geq s+1) \\
	&\geq 2 \quad (\text{since } s\geq 4) \\
	&> 0.
\end{align*}
Therefore, $g(n)$ is increasing on the interval $\left[\frac{2\delta s-3\delta+3}{6},+\infty\right)$. 
Note that
$\frac{2\delta s-3\delta+3}{6} \leq \frac{2\delta^{2}-5\delta+3}{6} < \frac{1}{3}\delta^{2}+3 \leq n$ by  $s\leq\delta-1.$ 
We conclude that
\begin{align*}
	g(n) &\geq g\left(\frac{1}{3}\delta^{2}+3\right) \\
	&= \left(\frac{s}{3}-\frac{2}{3}\right)\delta^{4} + \left(-\frac{2}{3}s^{2}+\frac{7}{3}s-\frac{4}{3}\right)\delta^{3} + \left(\frac{2}{3}s^{3}-\frac{17}{3}s^{2}+13s-10\right)\delta^{2} \\
	&\quad + (3s^{3}-14s^{2}+25s-13)\delta - s^{4}+7s^{3}-21s^{2}+33s-24 \\
	&\geq \frac{1}{3}s^{5}-\frac{4}{3}s^{4}+\frac{4}{3}s^{3}+34s-49 \quad (\text{since } \delta\geq s+1) \\
	&\geq 62 \quad (\text{since } s\geq 4) \\
	&> 0.\tag{23}
\end{align*}
From (21)$\sim(23)$, we get $\varphi_2(x) >0$ when $s \geq 4$.\\
For \( s = 3 \), we can directly get that
\begin{align*}
\varphi_2(x) &= (3n - 5\delta + 1)x + \delta n - n - 2\delta^2 + 4\delta - 3 \\
&\geq 3n^2 - (\delta + 9)n - 7\delta^2 + 20\delta - 6 \quad (\text{since } x \geq n + \delta - 3) \\
&\geq 177\delta^2 - 381\delta + 204 \quad (\text{since } n \geq 8\delta - 7) \\
&\geq 150 \quad (\text{since } \delta \geq 2) \\
&> 0.
\end{align*}\par
Togehter with $s \leq \delta - 1$ and $\varphi_2(x) >0$, we infer that $\varphi(B_{\Pi}^{3},x) <\varphi(B_{\Pi}^*,x)$ for $x \geq n+\delta-3$.
Notice that $\rho_{\mathcal{D}}(G_{*})\geq n+\delta-3$ by (16), we obtain that $\rho_{\mathcal{D}}(G_{3}) \geq \rho_{\mathcal{D}}(G_{*})$.
From the above disussion, we have
\begin{center}
	$\rho_{\mathcal{D}}(G) \geq \rho_{\mathcal{D}}(G_{1}) 
	\geq \rho_{\mathcal{D}}(G_{3}) 
	> \rho_{\mathcal{D}}(G_{*}) = \rho_{\mathcal{D}}(K_{\delta} \vee (K_{n-2\delta+1} \cup (\delta-1)K_{1})),$
\end{center}
which contradicts with $\rho_\mathcal{D}(G) \leq \rho_\mathcal{D}(K_\delta \vee (K_{n-2\delta+1} \cup (\delta - 1)K_1))$.\\
\textbf{Case 4.} $s=2$\par
In this case, $G_{3}=K_{2}\vee(K_{n-\delta-1}\cup K_{\delta-1})$. 
Let $\mathbf{x}$ be the Perron eigenvector of $G_{*}$. 
That is,
$\mathcal{D}(G_{*})\mathbf{x}=\rho_{\mathcal{D}}(G_{*})\mathbf{x}$. 
Since $\mathbf{x}$ is constant on each part corresponding to an equitable partition, let
\[
\mathbf{x}=(\underbrace{a,\cdots,a}_{\delta-1},\underbrace{b,\cdots,b}_{\delta},\underbrace{c,\cdots,c}_{n-2\delta+1})^{T},
\]
where $a,b,c\in\mathbb{R}^{+}$. 
In order to compare the appropriate spectral radii, we refine the partition, and write $\mathbf{x}$ as follows:
\[
\mathbf{x}=(\underbrace{a,\cdots,a}_{\delta-1},\underbrace{b,b}_{2},\underbrace{b,\cdots,b}_{\delta-2},\underbrace{c,\cdots,c}_{n-2\delta+1})^{T}.
\]
Thus, $\mathcal{D}(G_{3})-\mathcal{D}(G_{*})$ is partitioned as
\[
\bordermatrix{
	& \delta-1 & 2 & \delta-2 & n-2\delta+1 \cr
	\delta-1 & I-J & 0 & J & 0 \cr
	2 & 0 & 0 & 0 & 0 \cr
	\delta-2 & J & 0 & 0 & 0 \cr
	n-2\delta+1 & 0 & 0 & 0 & 0
},
\]
where $J$ is an all-ones matrix of appropriate size and $I$ is an identity matrix. 
Then we have
\begin{align*}
	\rho_{\mathcal{D}}(G_{3})-\rho_{\mathcal{D}}(G_{*}) &\geq \mathbf{x}^{T}(\mathcal{D}(G_{3})-\mathcal{D}(G_{*}))\mathbf{x} \\
	&= (\delta-1)(\delta-2)a(2b-a). \tag{24}
\end{align*}
Note that $\delta\geq s+1=3$ and $a>0$, it suffices to show that $2b-a>0$. \par
Recall that
\[
B^{*}_{\Pi}=\begin{bmatrix}
	n-2\delta & \delta & 2(\delta-1) \\
	n-2\delta+1 & \delta-1 & \delta-1 \\
	2(n-2\delta+1) & \delta & 2(\delta-2)
\end{bmatrix}
\]
is the equitable quotient matrix of $\mathcal{D}(G_{*})$. 
According to Lemma 2.5, we get $B^{*}_{\Pi}\mathbf{x}=\rho_{\mathcal{D}}(G_{*})\mathbf{x}$.
Hence
\[
\begin{cases}
	\rho_{\mathcal{D}}(G_{*})a = (n-2\delta)a + \delta b + 2(\delta-1)c & \text{(a)} \\
	\rho_{\mathcal{D}}(G_{*})b = (n-2\delta+1)a + (\delta-1)b + (\delta-1)c & \text{(b)} \\
	\rho_{\mathcal{D}}(G_{*})c = 2(n-2\delta+1)a + \delta b + 2(\delta-2)c & 
\end{cases} .\tag{25}
\]
From $2\text{(b)}-\text{(a)}$ of (25), we have $b = \dfrac{\rho_{\mathcal{D}}(G_{*})+n-2\delta+2}{2\rho_{\mathcal{D}}(G_{*})+\delta-2}a$. 
Thus
\begin{align*}
	2b - a &= \dfrac{2(\rho_{\mathcal{D}}(G_{*})+n-2\delta+2)}{2\rho_{\mathcal{D}}(G_{*})+\delta-2}a - a \\
	&= \dfrac{2n-5\delta+6}{2\rho_{\mathcal{D}}(G_{*})+\delta-2}a \\
	&\geq \dfrac{11\delta-8}{2\rho_{\mathcal{D}}(G_{*})+\delta-2}a \quad (\text{since } n\geq 8\delta-7).\tag{26}
\end{align*}
Note that $\rho_{\mathcal{D}}(G_{*}) > n+\delta-3 \geq 9\delta-10 > 0$ and $\delta\geq 2$, we can directly get that $2b-a>0$ from (26).\par
Concluding the above results, we conclude
\begin{center}
	$\rho_{\mathcal{D}}(G) \geq \rho_{\mathcal{D}}(G_{1}) 
	\geq \rho_{\mathcal{D}}(G_{3}) 
	> \rho_{\mathcal{D}}(G_{*})= \rho_{\mathcal{D}}(K_{\delta} \vee (K_{n-2\delta+1} \cup (\delta-1)K_{1})),$
\end{center}
a contradiction.\par
This completes the proof of Theorem 1.2.

\section*{Declaration of competing interest}

The authors declare that they have no known competing financial interests or personal relationships that could have appeared to influence the work reported in this paper.

\section*{Data availability}

No data was used for the research described in the article.

\section*{Acknowledgments}

We thank the anonymous reviewers for careful reading and their helpful comments.


\begin{thebibliography}{99}\addtolength{\itemsep}{-1.1ex}
\bibitem{01}M. Aouchiche, P. Hansen, Distance spectra of graphs: a survey, Linear Algebra Appl. 458 (2014) 301-386.
\bibitem{02} X. Chen, Q. Ji, M. Liu, Reducing Vizing's 2-factor conjecture to Meredith extension of critical graphs, Graphs Combin. 36 (2020) 1585-1591.
\bibitem{03}H. Fleischner, Spanning eulerian subgraphs, the Splitting Lemma, and Petersen's Theorem, Discrete Math. 101 (1992) 33-37. 
\bibitem{04} S. Fujita, K. Takazawa, The independent even factor problem, SIAM J. Discrete Math. 22 (2008) 1411-1427.
\bibitem{05} C. Godsil, Algebraic Combinatorics, Chapman and Hall Mathematics Series, New York, 1993.
\bibitem{06} C. Godsil, G. Royle, Algebraic graph theory, in: Graduate Texts in Mathematics, Vol.207, Springer-Verlag, New York, 2011.
\bibitem{07} R. Horn, C. Johnson, Matrix Analysis, Cambridge University Press, Cambridge (1985).
\bibitem{08} T. Király, M. Makai, On polyhedra related to even factors, in: G. Nemhauser, D. Bienstock (Eds.), Integer Programming and Combinatorial Optimization, Proc. 10th IPCO, in: Lecture Notes in Comput. Sci., vol. 3064, Springer-Verlag, 2004, pp. 416-430.
\bibitem{09}Y. Kobayashi, K. Takazawa, Even factors, jump systems, and discrete convexity, J. Combin. Theory Ser. B 99 (2009) 139-161.
\bibitem{10} S. Lv, L. Xiong, Even factors with a bounded number of components in iterated line graphs, Sci. China Math. 60 (2017) 177-188.
\bibitem{11}G. Pap, L. Sógŏ, On the maximum even factor in weakly symmetric graphs, J. Combin. Theory Ser. B 90 (2004) 201-213.
\bibitem{12} E. Steffen, I. Wolf, Even factors in edge-chromatic-critical graphs with a small number of divalent vertices, Graphs Combin. 38 (2022) 104.
\bibitem{13} Y. Shen, L. You, M. Zhang, S. Li, On a conjecture for the signless Laplacian spectral radius of cacti with given matching number, Linear Multilinear Algebra 65 (2017) 457-474.
\bibitem{14}L. Xiong, Characterization of forbidden subgraphs for the existence of even factors in a graph, Discrete Appl. Math. 223 (2017) 135-139.
\bibitem{15} L. Xiong, The existence of even factors in iterated line graphs, Discrete Appl. Math. 308 (2023) 591-594.
\bibitem{16}Z. Yuan, M. Kano, Strong Tutte type conditions and factors of graphs, Discuss. Math. Graph Theory 40 (2020) 1057-1065.
\bibitem{17} L. You, M. Yang, W. So, W. Xi, On the spectrum of an equitable quotient matrix and its application, Linear Algebra Appl. 577 (2019) 21-40.
\bibitem{18}S. Zhou, Q. Bian, J. Wu, Sufficient conditions for even factors in graphs, arXiv preprint arXiv:2510.10600 (2025).
\bibitem{19}L. Zheng, S. Li, X. Luo, G. Wang, Some sufficient conditions for a graph with minimum degree to be k-factor-critical, Discrete Appl. Math. 348 (2024) 279-291.
\bibitem{20}L. Zhang, L. Xiong, Characterizing forbidden pairs for the existence of even factors, Discrete Math. 348 (2025) 114384.
\bibitem{21}S. Zhou, L. Zhang, Signless Laplacian spectral radius for k-extendable graphs, Filomat 2 (2014) 659-667.
	
    

\end{thebibliography}
\end{document}